\documentclass{article}
\usepackage{mathrsfs}
\usepackage{bbm}
\usepackage{amssymb,amsmath,amsthm,cases}
\usepackage[noadjust]{cite}
\allowdisplaybreaks
\begin{document}
\title {On the octonionic Bergman kernel}
\author{Jinxun Wang\thanks{\scriptsize Department of
Applied Mathematics, Guangdong University of Foreign Studies, Guangzhou 510006,
China. E-mail: wjx@gdufs.edu.cn}\,, ~Xingmin
Li\thanks{\scriptsize School of Computer Sciences, South China
Normal University, Guangzhou 510631, China. E-mail: lxmin57@163.com}}
\date{}
\maketitle
\newtheorem{defi}{Definition}[section]
\newtheorem{theo}{Theorem}[section]
\newtheorem{prop}{Proposition}[section]
\newtheorem{lemm}{Lemma}[section]
\newtheorem{coro}{Corollary}[section]

\noindent\textbf{Abstract:} By introducing a suitable new definition for
the inner product on the octonionic Bergman space, we determine the explicit
form of the octonionic Bergman kernel, in the framework of octonionic analysis
which is non-commutative and non-associative.

\section{Introduction}
In complex analysis the Szeg\"{o} kernel and Bergman kernel are
well-known, which had also been generalized into Clifford analysis
(including quaternionic analysis as a special case, see \cite{BDS}).
But in octonionic analysis the existence of such kernels is still
unknown, let alone the explicit expressions. The difficulty arises
mainly because the octonion algebra (Cayley algebra) is non-associative.

The motivation for us to consider this kind of problem is that we
want to unify the formulation of the analytic function theory in the
largest normed division algebra over $\mathbb{R}$, namely, in octonions
$\mathbb{O}$ (including complex numbers,
quaternions as its special cases).

Recall that in complex analysis the Bergman space on the unit disc
is defined to be the collection of functions that are holomorphic
and square integrable on the unit disc. This definition can be
naturally generalized to octonionic analytic functions. Since the
Cayley algebra is non-commutative, there exist two different but
symmetric octonionic analytic function theory. In this paper we
focus on the left octonionic analytic functions, and we denote by
$\mathcal{B}^2(B)$ the corresponding octonionic Bergman space,
where $B$ is the unit ball in $\mathbb{R}^8$ centered at origin. A nature problem
comes: Does the octonionic Bergman kernel exist? and what is it? Of
course this problem is closely related to the definition of the
associated inner product. Usually the inner product of two Bergman
functions $f$ and $g$ is defined to be the integral of
$f\overline{g}$ on $B$. Since the octonions is non-associative, the
usual definition is no longer valid to guarantee the existence of
the kernel. We thus need to give a new definition.

\begin{defi}[inner product on $\mathcal{B}^2(B)$]
Let $f, g\in\mathcal{B}^2(B)$, we define
$$(f,g)_B:=\frac{1}{\omega_8}\int_B\left(\overline{g(x)}\frac{\overline{x}}{|x|}\right)
\left(\frac{x}{|x|}f(x)\right)dV,$$ where $\omega_8$ is the surface
area of the unit sphere in $\mathbb{R}^8$, $dV$ is the volume
element on $B$.
\end{defi}

Note that this modified inner product is real-linear and conjugate
symmetric. The induced norm
$$\|f\|_B^2:=(f,f)_B=\frac{1}{\omega_8}\int_B|f|^2dV$$
coincides with the norm induced by the usual inner product.

We can now state the main theorem of this paper.
\begin{theo}\label{main}
Let$$B(x,a)=\frac{\left(6(1-|a|^2|x|^2)+2(1-\overline{x}a)\right)(1-\overline{x}a)}
{|1-\overline{x}a|^{10}},$$ then $B(\cdot,a)$ is the desired
octonionic Bergman kernel, i.e., $B(\cdot,a)\in\mathcal{B}^2(B)$,
and for any $f\in\mathcal{B}^2(B)$ and any $a\in B$, there holds the
following reproducing formula
$$f(a)=(f,B(\cdot,a))_B.$$
\end{theo}

The rest of the paper is organized as follows. In Section 2 we give a brief review on
the octonion algebra and octonionic analysis. In Section 3 we will exploit our new idea in defining the
structure of the inner product to investigate the octonionic Szeg\"{o} kernel for the unit
ball in $\mathbb{R}^8$. Section 4 is then devoted to the proof of our main result Theorem \ref{main}. In the
last section we point out that the Bergman kernel can be unified in one form in both complex analysis and
hyper-complex analysis.

\section{The octonions and the octonionic analysis}
\subsection{The octonions}
If an algebra $\mathbb{A}$ is meanwhile a normed vector space, and
its norm ``$\|\cdot\|$'' satisfies $\|ab\|=\|a\|\|b\|$, then we call
$\mathbb{A}$ a normed algebra. If $ab=0$ ($a, b\in \mathbb{A}$)
implies $a=0$ or $b=0$, then we call $\mathbb{A}$ a division
algebra. Early in 1898, Hurwitz had proved that the real numbers
$\mathbb{R},$ complex numbers $\mathbb{C},$ quaternions $\mathbb{H}$
and octonions $\mathbb{O}$ are the only normed division algebras
over $\mathbb{R}$ (\cite{Hur}), with the imbedding relation
$\mathbb{R}\subseteq \mathbb{C}\subseteq \mathbb{H}\subseteq
\mathbb{O}$.

As the largest normed division algebra, octonions, which are also called
Cayley numbers or the Cayley algebra, were discovered by
John T. Graves in 1843, and then by Arthur Cayley in 1845
independently. Octonions are an 8 dimensional algebra over
$\mathbb{R}$ with the basis $e_0,e_1,\ldots,e_7$ satisfying
$$e_0^2=e_0,~e_ie_0=e_0e_i=e_i,~e_i^2=-1,~\text{for}~i=1,2,\ldots,7.$$
So $e_0$ is the unit element and can be identified with $1$.
Denote
$$W=\{(1,2,3),(1,4,5),(1,7,6),(2,4,6),(2,5,7),(3,4,7),(3,6,5)\}.$$
For any triple $(\alpha,\beta,\gamma)\in W$, we set $$e_\alpha
e_\beta=e_\gamma=-e_\beta e_\alpha,\quad e_\beta
e_\gamma=e_\alpha=-e_\gamma e_\beta,\quad e_\gamma
e_\alpha=e_\beta=- e_\alpha e_\gamma.$$
Then by distributivity for any $x=\sum_0^7 x_ie_i$,
$y=\sum_0^7 y_je_j \in \mathbb{O}$, the multiplication $xy$ is defined to be
$$xy:=\sum_{i=0}^7\sum_{j=0}^7x_iy_je_ie_j.$$

For any $x=\sum_0^7 x_ie_i \in \mathbb{O}$, $\mbox{Re}\,x:=x_{0}$ is
called the scalar (or real) part of $x$ and
$\overrightarrow{x}:=x-\mbox{Re}\,x$ is called its vector part.
$\overline{x}:=\sum_0^7x_i\overline{e_i}=x_0-\overrightarrow{x}$
and $|x|:=(\sum_0^7x_i^2)^\frac{1}{2}$ are respectively the conjugate and
norm (or modulus) of $x$, they satisfy: $|xy|=|x||y|,$
$x\overline{x}=\overline{x}x=|x|^2,$
$\overline{xy}=\overline{y}\,\overline{x}$ $(x,y\in \mathbb{O}).$ So if
$x\neq0,$ $x^{-1}=\overline{x}/{|x|^2}$ gives the inverse of $x.$

Octonionic multiplication is neither commutative nor associative. But
the subalgebra generated by any two elements is associative, namely, The octonions
are alternative. $[x, y, z]:=(xy)z-x(yz)$ is called the associator of $x, y, z\in
\mathbb{O},$ it satisfies (\cite{B1, Jaco})
$$
[x,y,z]=[y,z,x]=-[y,x,z], \quad [x,x,y]=[\overline{x},x,y]=0.
$$

\subsection{The octonionic analysis}
As a generalization of complex analysis and quaternionic
analysis to higher dimensions, the study of octonionic analysis
was originated by Dentoni and Sce in 1973 (\cite{DS}), and it was not until
1995 that it began to be systematically investigated by Li et al (\cite{Li}).
Octonionic analysis is a function theory on
octonionic analytic (abbr. $\mathbb{O}$-analytic) functions. Suppose $\Omega$ is
an open subset of $\mathbb{R}^8$,
$f=\sum_0^7f_je_j\in C^1(\Omega,\mathbb{O})$ is an octonion-valued function, if
$$Df=\sum_{i=0}^{7}e_{i}\frac{\partial f}{\partial x_{i}}
=\sum_{i=0}^{7}\sum_{j=0}^{7}\frac{\partial f_j}{\partial x_{i}}e_ie_j=0$$
$$\left(fD=\sum_{i=0}^{7} \frac{\partial f}{\partial x_{i}}e_{i}=
\sum_{i=0}^{7}\sum_{j=0}^{7}\frac{\partial f_j}{\partial x_{i}}e_je_i=0\right),$$ then $f$
is said to be left (right) $\mathbb{O}$-analytic in $\Omega$, where
the generalized Cauchy--Riemann operator $D$ and its conjugate $\overline{D}$ are defined
by $$D:=\sum_{i=0}^7e_i\frac{\partial}{\partial x_i},~~\overline{D}:=
\sum_{i=0}^7\overline{e_i}\frac{\partial}{\partial x_i}$$ respectively.
A function $f$ is $\mathbb{O}$-analytic means that $f$ is meanwhile left
$\mathbb{O}$-analytic and right $\mathbb{O}$-analytic. From
$$\overline{D}(Df)=(\overline{D}D)f=\triangle
f=f(D\overline{D})=(fD)\overline{D},$$ we know that any left (right)
$\mathbb{O}$-analytic function is always harmonic. In the sequel,
unless otherwise specified, we just consider the left
$\mathbb{O}$-analytic case as the right $\mathbb{O}$-analytic case is essentially the same.
A Cauchy-type integral formula and a Laurent-type series for this setting are:

\begin{lemm}[Cauchy's integral formula, see \cite{DS,LP2}]
Let $\mathcal{M}\subset\Omega$ be an 8-dimensional, compact differentiable and oriented
manifold with boundary. If $f$ is left $\mathbb{O}$-analytic in $\Omega$, then
$$f(x)=\frac{1}{\omega_8}\int_{y\in\partial\mathcal{M}}E(y-x)(d\sigma_yf(y)),\quad x\in
\mathcal{M}^o,$$
where $E(x)=\frac{\overline{x}}{|x|^8}$ is the octonionic Cauchy kernel, $d\sigma_y
=n(y)dS$, $n(y)$ and $dS$ are respectively the outward-pointing unit normal vector
and surface area element on $\partial\mathcal{M}$,
$\mathcal{M}^o$ is the interior of $\mathcal{M}$.
\end{lemm}

\begin{lemm}[Laurent expansion, see \cite{LLW,LZP2}]
Let $\mathcal{D}$ be an annular domain in $\mathbb{R}^8$.
If $f$ is left $\mathbb{O}$-analytic in $\mathcal{D}$, then
$$f(x)=\sum_{k=0}^\infty P_kf(x)+\sum_{k=0}^\infty Q_kf(x),
\quad x\in\mathcal{D},$$ where $P_kf$ and $Q_kf$ are respectively the inner and outer
spherical octonionic-analytics of order $k$ associated to $f$.
\end{lemm}

Octonionic analytic functions have a close relationship with the
Stein--Weiss conjugate harmonic systems. If the
components of $F$ consist a Stein--Weiss conjugate
harmonic system on $\Omega\subset\mathbb{R}^8$,
then $\overline{F}$ is $\mathbb{O}$-analytic on $\Omega$.
But conversely this is not true (\cite{LP1}).
For more information and recent progress about octonionic analysis, we
refer the reader to \cite{Li, LPQ1, LPQ2, LW, LZP1, LLW, WLL}.

\section{The octonionic Szeg\"{o} kernel}
To see how our new definition works, let us check the octonionic
Szeg\"{o} kernel for the unit ball in $\mathbb{R}^8$.

Recall that on the unit ball the octonionic Hardy space
$\mathcal{H}^2(B)$ consists of the left octonionic analytic
functions whose mean square value on the sphere is bounded for
radius $r\in[0,1)$. For any $f\in\mathcal{H}^2(B)$, according to the
Cauchy's integral formula, for all $a\in B$ there holds
\begin{align*}
f(a)&=\frac{1}{\omega_8}\int_{x\in
S^7}\frac{\overline{x}-\overline{a}}{|x-a|^8}(xf(x))dS
\\&=\frac{1}{\omega_8}\int_{x\in S^7}\left(
\frac{\overline{1-\overline{x}a}}{|1-\overline{x}a|^8}\overline{x}
\right)(xf(x))dS,
\end{align*}
where $S^7=\partial B$ is the unit sphere, $dS$ is the area element
on $S^7$. If we define the inner product for $\mathcal{H}^2(B)$ to
be
\begin{align}\label{inner1}
(f, g)_{S^7}:=\frac{1}{\omega_8}\int_{S^7}(\overline{\eta
g(\eta)})(\eta f(\eta))dS=
\frac{1}{\omega_8}\int_{S^7}(\overline{g(\eta)}\overline{\eta})(\eta
f(\eta))dS,
\end{align}
and let
$$S(x,a)=\frac{1-\overline{x}a}{|1-\overline{x}a|^8},$$
then $S(\cdot,a)\in\mathcal{H}^2(B)$, and the Cauchy's integral
formula can be rewritten as
$$f(a)=(f,S(\cdot,a))_{S^7}.$$
We call $S(\cdot,a)$ the octonionic Szeg\"{o} kernel.

Denote by $L^2(S^7)$ the space of square integrable
(octonion-valued) functions on the unit sphere, for which we define
its inner product to be the same as that in (\ref{inner1}). We have
\begin{prop}\label{propo}
Let $f, g\in L^2(S^7)$ be associated with the spherical
octonionic-analytics expansions:
$$f(\omega)=\sum_{k=0}^\infty(P_kf(\omega)+Q_kf(\omega)),\quad
g(\omega)=\sum_{k=0}^\infty(P_kg(\omega)+Q_kg(\omega)),\quad\omega\in
S^7.$$ Then
\begin{align*}
(f, g)_{S^7}=&\sum_{k=0}^\infty\left((P_kf, P_kg)_{S^7}+(Q_kf,
Q_kg)_{S^7}\right)\\&+\sum_{k=0}^\infty\left((P_kf,
Q_{k+1}g)_{S^7}+(Q_{k+1}f, P_kg)_{S^7}\right).
\end{align*}
\begin{proof}
From $$\triangle(xP_kf(x))=x\triangle(P_kf(x))+2D(P_kf(x))=0,$$ we
can easily see that the restriction of $xP_kf(x)$ on
$S^7$ is a spherical harmonic of order $k+1$. Similarly, the restriction of
$xQ_kf(x)$ on $S^7$ is a spherical harmonic of order $k$.
The proposition immediately follows by the fact that spherical
harmonics of different orders are mutually orthogonal.
\end{proof}
\end{prop}

Thus we get
\begin{coro}
Let $f\in L^2(S^7)$ be associated with the spherical
octonionic-analytics expansion
$$f(\omega)=\sum_{k=0}^\infty(P_kf(\omega)+Q_kf(\omega)),\quad
\omega\in S^7.$$ Then
$$
\|f\|^2_{S^7}=\sum_{k=0}^\infty\left(\|P_kf\|^2_{S^7}+\|Q_kf\|^2_{S^7}\right)
+\sum_{k=0}^\infty2{\rm Re}\left((P_kf, Q_{k+1}f)_{S^7}\right).
$$
\end{coro}

\noindent\textbf{Remark:} Proposition \ref{propo} is similar to the Parseval's
theorem. It is worthwhile to note that this version is a bit different
from that in Clifford analysis where the second part in the summation
vanishes (\cite{BDS}), here $(P_kf, Q_{k+1}g)_{S^7}$ may not be zero.
Below we give a counter-example. Let
$$f(x)=x_1-x_0e_1,$$
$$g(x)=\frac{\overline{x}}{|x|^{12}}(x_1x_2e_4+x_0x_2e_5+x_0x_1e_6).$$
Then $P_1f=f$, $Q_2g=g$, but
$$(P_1f, Q_2g)_{S^7}=\frac{-2e_6}{\omega_8}\int_{S^7}x_0^2x_1^2dS\neq
0.$$

\section{Derivation of the octonionic Bergman kernel}
In this section we will prove Theorem \ref{main}. For the main idea
we use in the proof one can also refer to \cite{BDS}.

\begin{proof}[Proof of Theorem \ref{main}]
By definition it is straightforward that
$$(f,g)_B=\int_0^1r^7(f_r,g_r)_{S^7}dr,$$
where $f_r(\eta)=f(r\eta)$, $\eta\in S^7$. Together with Proposition
\ref{propo}, we get
\begin{align*}
(f,g)_B&=\sum_{k=0}^\infty(P_kf,P_kg)_B\\
&=\sum_{k=0}^\infty\int_0^1r^{2k+7}(P_kf,P_kg)_{S^7}dr\\
&=\sum_{k=0}^\infty(2k+8)^{-1}(P_kf,P_kg)_{S^7}.
\end{align*}
Therefore, $f\in\mathcal{B}^2(B)$ if and only if $f$ is left
octonionic analytic in $B$ and
$$\|f\|_B^2=\sum_{k=0}^\infty(2k+8)^{-1}\|P_kf\|_{S^7}^2<\infty.$$
From this viewpoint, if $f\in\mathcal{H}^2(B)$, then
$$\sqrt{T}f:=\sum_{k=0}^\infty\sqrt{2k+8}P_kf\in\mathcal{B}^2(B).$$
Similarly, if $g$ is left octonionic analytic in $B_R$ (the ball
centered at the origin of radius $R$, with $R>1$), then
$\sqrt{T}g\in\mathcal{B}^2(B_{R'})$, with $1\leq R'<R$. Consequently,
$$Tg:=\sqrt{T}^2g=\sum_{k=0}^\infty(2k+8)P_kg
\in\mathcal{B}^2(B_{R'}),~1\leq R'<R.$$

Now, assume $f\in\mathcal{B}^2(B)$, when $|a|<r$ we have
\begin{align}
f(a)&=\frac{1}{\omega_8}\int_{\partial B_r}
\frac{\overline{x}-\overline{a}}{|x-a|^8}d\mu_xf(x)\nonumber\\
&=\frac{r^7}{\omega_8}\int_{S^7}\frac{r\overline{\eta}-\overline{a}}
{|r\eta-a|^8}(\eta f(r\eta))dS\nonumber\\
&=\lim_{r\rightarrow 1^-}\frac{r^7}{\omega_8}\int_{S^7}
\frac{r\overline{\eta}-\overline{a}}
{|r\eta-a|^8}(\eta f(r\eta))dS\nonumber\\
&=\lim_{r\rightarrow 1^-}r^7(f_r, S^r(\cdot,a))_{S^7}\label{proof1},
\end{align}
where $$S^r(x,a)=\frac{r-\overline{x}a}{|r-\overline{x}a|^8}.$$
Since $S^r(x,a)$ is left octonionic analytic in $B_{r/|a|}$
($r/|a|>1$) with respect to $x$, we have
$$TS^r(\cdot,a)=\sum_{k=0}^\infty(2k+8)P_kS^r(\cdot,a)\in\mathcal{B}^2(B).$$
So,
\begin{align}
(f_r,TS^r(\cdot,a))_B&=\sum_{k=0}^\infty(2k+8)^{-1}
(P_kf_r,(2k+8)P_kS^r(\cdot,a))_{S^7}\nonumber\\&=\sum_{k=0}^\infty
(P_kf_r,P_kS^r(\cdot,a))_{S^7}\nonumber\\&=(f_r,S^r(\cdot,a))_{S^7}\label{proof2}.
\end{align}
By (\ref{proof1}) and (\ref{proof2}) we get
$$f(a)=\lim_{r\rightarrow 1^-}r^7(f_r,TS^r(\cdot,a))_B
=(f,TS(\cdot,a))_B,$$ where $S(\cdot,a)$ is the octonionic Szeg\"{o}
kernel. We can now see that the octonionic Bergman kernel $B(x,a)$
is
$$B(x,a)=TS(x,a)=\sum_{k=0}^\infty(2k+8)P_kS(x,a).$$

The remaining thing we need to do is to evaluate the above
summation. To this end, first note that
$$S(x,a)=\mathcal{K}(E(x,\overline{a})),$$
where $E(x,a)=\frac{\overline{x}-\overline{a}}{|x-a|^8}$ ($|x|>1$),
$\mathcal{K}f:=E(x,0)f(x^{-1})$ is the Kelvin inversion. So,
$$P_kS(x,a)=\mathcal{K}(Q_kE(x,\overline{a}))=\overline{x}
\overline{Q_kE(x,a)}|x|^{2k+6}.$$ Define the adjoint operator $A$ as
follows:
$$(Af)(x):=\overline{D}(|x|^{-6}\overline{f}(x/|x|^2)),$$
then it is easy to show that
$$A(Q_kE(x,a))=(2k+8)\overline{x}
\overline{Q_kE(x,a)}|x|^{2k+6}.$$
Hence,
\begin{align*}
B(x,a)&=\sum_{k=0}^\infty A(Q_kE(x,a))\\
&=A\left(\sum_{k=0}^\infty Q_kE(x,a)\right)\\
&=A(E(x,a))
\\&=\overline{D}_x\left(\frac{x-a|x|^2}{|1-\overline{x}a|^8}\right)
\\&=\frac{\left(6(1-|a|^2|x|^2)+2(1-\overline{x}a)\right)(1-\overline{x}a)}
{|1-\overline{x}a|^{10}}.
\end{align*}
The proof of Theorem \ref{main} is complete.
\end{proof}

\section{Final remarks}
By direct computation one can show that
$$\overline{B(x,a)}\overline{x}=
\overline{D}_a\left(\frac{1-|a|^2|x|^2}{|1-a\overline{x}|^8}\right).$$
In fact, similar formulas also hold in both complex analysis and
Clifford analysis. We therefore can unify the reproducing formulas
in complex and hyper-complex contexts. Let $\mathscr{A}$ denote the
complex algebra or hyper-complex algebra, i.e., $\mathscr{A}$ may
refer to complex numbers $\mathbb{C}$, quaternions $\mathbb{H}$,
octonions $\mathbb{O}$, or Clifford algebra $\mathscr{C}$. Assume
that the dimension of $\mathscr{A}$ is $m$. Then for any function
$f$ which belongs to the Bergman space $\mathcal{B}^2(B_m)$ and any
point $a\in B_m$ ($B_m$ is the unit ball centered at origin
in $\mathbb{R}^m$), there holds
\begin{align*}
f(a)&=(f,B(\cdot,a))_{B_m} \\
&=\frac{1}{\omega_m}\int_{B_m}\left(\overline{B(x,a)}\frac{\overline{x}}{|x|}\right)
\left(\frac{x}{|x|}f(x)\right)dV
\\&=\frac{1}{\omega_m}\int_{B_m}
\overline{D}_a\frac{1-|a|^2|x|^2}{|1-a\overline{x}|^m}
\left(\frac{x}{|x|^2}f(x)\right)dV
\\&=\frac{1}{\omega_m}\int_{B_m}
\frac{\left((m-2)(1-|a|^2|x|^2)+2(1-\overline{a}x)\right)(\overline{x}-|x|^2\overline{a})}
{|1-\overline{x}a|^{m+2}}\left(\frac{x}{|x|^2}f(x)\right)dV,
\end{align*}
where $\omega_m$ is the surface
area of the unit sphere in $\mathbb{R}^m$, $dV$ is the volume
element on $B_m$, and $D$ is the generalized Cauchy--Riemann operator in the respective
context.

\vskip 0.8cm \noindent{\Large\textbf{Acknowledgements}}

\vskip 0.3cm \noindent
This work was supported by the Scientific
Research Grant of Guangdong University of Foreign Studies for Introduction of Talents
(No. 299--X5122145), the Research Grant of Guangdong University of Foreign
Studies for Young Scholars (No. 299--X5122199), and the Foundation for
Young Innovative Talents in Higher Education of Guangdong, China (No. 2015KQNCX037).

\end{document}